\documentclass[11pt]{article}
\usepackage{latexsym,amsmath,color,amsthm,amssymb,epsfig,graphicx,mathrsfs}
\usepackage{fullpage,sectsty}
\usepackage[titles]{tocloft}
\usepackage[compact]{titlesec}
\titleformat{\subsection}[runin]{\normalfont\bfseries}{\thesubsection}{0.5em}{}
\titlespacing{\section}{0pt}{*4}{*2}
\titlespacing{\subsection}{0pt}{*3}{*1}

\def\qed{\hfill \ifhmode\unskip\nobreak\fi\quad\ifmmode\Box\else$\Box$\fi\\ }
\def\ex{{\rm ex}}

\newtheorem{thrm}{Theorem}[section]
\newtheorem{lemm}[thrm]{Lemma}
\newtheorem{propo}[thrm]{Proposition}
\newtheorem{coro}[thrm]{Corollary}
\newtheorem{defi}[thrm]{Definition}
\newtheorem{conjecture}[thrm]{Conjecture}

\newcommand{\thm}{\begin{thrm}}
\newcommand{\xthm}{\end{thrm}}
\newcommand{\lem}{\begin{lemm}}
\newcommand{\xlem}{\end{lemm}}
\newcommand{\prf}{\begin{proof}}
\newcommand{\xprf}{\end{proof} \vspace{-0.3in}}
\newcommand{\prop}{\begin{propo}}
\newcommand{\xprop}{\end{propo}}
\newcommand{\cor}{\begin{coro}}
\newcommand{\xcor}{\end{coro}}
\newcommand{\defn}{\begin{defi}}
\newcommand{\xdefn}{\end{defi}}
\newcommand{\conj}{\begin{conjecture}}
\newcommand{\xconj}{\end{conjecture}}

\renewcommand{\phi}{\varphi}

\newcommand{\PP}{\mathcal{P}}
\newcommand{\CC}{\mathcal{C}}
\begin{document}

\title{Tur\'an Problems and Shadows I: Paths and Cycles}

\author{
{\large{Alexandr Kostochka}}\thanks{
\footnotesize {University of Illinois at Urbana--Champaign, Urbana, IL 61801
 and Sobolev Institute of Mathematics, Novosibirsk 630090, Russia. E-mail: \texttt {kostochk@math.uiuc.edu}.
 Research of this author
is supported in part by NSF grants DMS-0965587 and DMS-1266016
and by
grant 12-01-00631 of the Russian Foundation for Basic Research.
}}
\and
{\large{Dhruv Mubayi}}\thanks{
\footnotesize {Department of Mathematics, Statistics, and Computer Science, University of Illinois at Chicago, Chicago, IL 60607.
E-mail:  \texttt{mubayi@uic.edu.}
Research partially supported by NSF grants DMS-0969092 and DMS-1300138.}}
\and{\large{Jacques Verstra\"ete}}\thanks{Department of Mathematics, University of California at San Diego, 9500
Gilman Drive, La Jolla, California 92093-0112, USA. E-mail: {\tt jverstra@math.ucsd.edu.} Research supported by NSF Grant DMS-1101489. }}
\maketitle

\begin{abstract}

A {\em $k$-path} is a hypergraph $P_k = \{e_1,e_2,\dots,e_k\}$ such that
$|e_i \cap e_j| = 1$ if $|j - i| = 1$
 and $e_i \cap e_j = \emptyset$ otherwise. A {\em $k$-cycle} is a hypergraph $C_k = \{e_1,e_2,\dots,e_k\}$
obtained from a $(k-1)$-path  $\{e_1,e_2,\dots,e_{k-1}\}$ by adding an edge $e_k$ that shares
one vertex with $e_1$, another vertex with $e_{k-1}$ and is disjoint from the other edges.

Let $\ex_r(n,G)$ be the maximum number
of edges in an $r$-graph with $n$ vertices not containing a given $r$-graph $G$.
We prove that for fixed $r \geq 3,k \geq 4$ and $(k,r) \neq (4,3)$, for large enough $n$:
\[ \ex_r(n,P_k) = \ex_r(n,C_k) = \binom{n}{ r} - \binom{n - \lfloor \frac{k-1}{2}\rfloor}{r} +
  \left\{\begin{array}{ll}
0 & \mbox{ if }k\mbox{ is odd}\\
{n - \lfloor \frac{k-1}{2}\rfloor - 2 \choose r - 2} & \mbox{ if }k\mbox{ is even}
\end{array}\right.
\]
and we characterize all the extremal $r$-graphs. We also solve the case $(k,r) = (4,3)$, which needs
a special treatment. The case $k  = 3$
was settled by Frankl and F\"{u}redi. 

This work is the next step in a long line of research beginning with  conjectures
of Erd\H os and S\'os from the early 1970's.
In particular, we extend the work (and settle a 
 conjecture) of F\"uredi, Jiang and
Seiver who solved this problem for $P_k$ when $r \ge 4$ and of F\"uredi and
Jiang who solved it for $C_k$ when $r \ge 5$. They used the delta system method, while
we use a novel approach which involves random sampling from the shadow of an $r$-graph.
\end{abstract}


\section{Introduction}

An $r$-uniform hypergraph, or simply {\em $r$-graph}, is a family of $r$-element subsets of a finite set.
Given a set  $\mathcal{F}$ of $r$-graphs, an  $\mathcal{F}$-{\em free $r$-graph} is an $r$-graph containing none
of the members of  $\mathcal{F}$.
Let the {\em Tur\' an number of $\mathcal{F}$}, $\ex_r(n,\mathcal{F})$,
denote the maximum number of edges in an  $\mathcal{F}$-free $r$-graph on $n$ vertices. 
 When $\mathcal{F}=\{F\}$ we write $\ex_r(n,F)$.
An  $n$-vertex  $\mathcal{F}$-free $r$-graph $H$ is {\em extremal for} $\mathcal{F}$ if $|H|= \ex_r(n, \mathcal{F})$.
In this paper we promote the idea of determining $\ex_r(n, \mathcal{F})$ for certain classes $\mathcal{F}$ by
 randomly sampling from the shadow of an $\mathcal{F}$-free $r$-graph $H$ and using Hall-type combinatorial lemmas to
 determine the structure of the shadow and hence the structure of $H$.
 This paper focuses solely on paths and cycles. Our next paper will consider more general structures.

\medskip


\subsection{Definitions of paths and cycles.\mbox{ }} There are several natural generalizations to hypergraphs of paths and cycles in graphs.
A {\em Berge $k$-cycle} is a hypergraph consisting of $k$ distinct edges $e_0,\ldots,e_{k-1}$ such that there exist
$k$ distinct vertices $v_0,v_1\ldots,v_{k-1}$ with $v_i\in e_{i-1}\cap e_i$ for all $i= 0,1,\ldots,k-1$ (indices count modulo $k$).
Let $\mathcal{B}\mathcal{C}_k$ denote the family of all Berge $k$-cycles.
A {\em minimal $k$-cycle} is a Berge cycle $\{e_0,e_1,\dots,e_{k-1}\}$ such that
$e_i \cap e_j \neq \emptyset$ if and only if $|j - i| = 1$ or $\{i,j\} = \{0,k-1\}$, and 
 no
vertex belongs to all edges. Let $\mathcal{C}_k$ denote the family of minimal $k$-cycles.
Furthermore, a {\em linear $k$-cycle}  is the member $C_k \in \mathcal{C}_k$ such that $|e_i\cap e_{i+1}| = 1$
for all $i=0,1,\ldots,k-1$.

Every  Berge (respectively, minimal and linear) {\em $k$-path} is obtained from a Berge
(respectively, minimal and linear)  $(k+1)$-cycle
 by deleting one edge.
The family  of Berge (respectively, minimal) $k$-paths is denoted by ${\cal BP}_k$ (respectively, ${\cal P}_k$).
The linear $k$-path is denoted by $P_k$. The most restricted structures above are linear $k$-cycles and $k$-paths.
We will refer to these simply as  {\em $k$-cycles} and  {\em $k$-paths}. In this paper, we study the extremal functions
for $k$-paths and $k$-cycles and  minimal  $k$-paths and $k$-cycles.

\subsection{The extremal function for $k$-cycles and $k$-paths.\mbox{ }} The extremal problem for $P_k$ has been studied extensively.
In the case of graphs, the Erd\H{o}s-Gallai Theorem~\cite{EG} shows $\ex(n,P_k) \leq \frac{k-1}{2}n$ and this is tight whenever $k|n$.
Frankl~\cite{Frankl} solved the simplest case for $r$-graphs, namely $\ex_r(n, P_2)$,
answering a question of Erd\H os and S\'os. As far as exact results are concerned, it appears that even the next smallest
case $\ex_r(n, P_3)$ was not determined until very recently. F\"uredi, Jiang and Seiver~\cite{FJS} determined
$\ex_r(n, P_k)$ precisely for all $r \ge 4$, $k \ge 3$ and $n$ large while also characterizing the extremal examples.
They conjectured a similar result for $r=3$.
In this paper, we prove their conjecture and determine the extremal structures for large $n$.

The extremal problem for $r$-graphs for $C_3$ is also well-researched~\cite{CK,FF}, indeed, the case $r=2$
 is precisely Mantel's theorem from 1907.  Frankl and F\"{u}redi~\cite{FF} showed that the unique  extremal $r$-graph
 on $[n]$ not containing $C_3$ consists of all edges containing some $x \in [n]$, for large enough $n$.
 For $r=k=3$ the exact result was proved for all $n \ge 6$ by Cs\'ak\'any  and Kahn~\cite{CK}.
 More recently, F\"uredi and Jiang~\cite{FJ} determined the extremal function for $C_k$ for all $k \ge 3$,
 $r \ge 5$ and large $n$; their results substantially extend earlier results of Erd\H os and settled a conjecture
 of the last two authors for $r \ge 5$. They used the delta system method.

 Our main result extends the F\"uredi-Jiang Theorem to the case of $r=3,4$.
To describe the result, we need some notation.
Let $[n] := \{1,2,\dots,n\}$, and for $L\subset [n]$ let $S_{L}^r(n)$
denote the $r$-graph on $[n]$ consisting of all $r$-element subsets of $[n]$
intersecting $L$.

\thm \label{main}
Let $r \geq 3$, $k\geq 4$, and $\ell = \lfloor \frac{k-1}{2} \rfloor$. For sufficiently large $n$, 
\[
\ex_r(n,P_k) = {n \choose r} - {n - \ell \choose r} + \left\{\begin{array}{ll}
0 & \mbox{ if }k\mbox{ is odd}\\
{n - \ell - 2 \choose r - 2} & \mbox{ if }k\mbox{ is even}
\end{array}\right.
\]
with equality only for $S_L^r(n)$ if $k$ is odd and $S_L^r(n) \cup F$ where $F$ is extremal for $\{P_2,2P_1\}$ on $n - \ell$ vertices.
The same result holds for $k$-cycles except  the case $(k,r) = (4,3)$, in which case
\[ \ex_3(n,C_4) = {n \choose r} - {n - 1 \choose r} + \max\{n-3,4\lfloor \tfrac{n-1}{4}\rfloor\}\]
with equality only for 3-graphs of the form $S_L^3(n) \cup F$ where $F$ is extremal for $P_2$ on $n - 1$ vertices.
  \xthm

{\bf Remarks.}

(1) By the Erd\H os-Ko-Rado Theorem~\cite{EKR}, $\ex_r(n-\ell,\{P_2,2P_1\}) = \binom{n-\ell-2}{r-2}$ for sufficiently large $n$,
 and a result of Erd\H os and S\' os~(see~\cite{Frankl}) gives $\ex_3(n-1,P_2) = \max\{n-3,4\lfloor \tfrac{n-1}{4}\rfloor\}$. These results account
 for the lower order terms in the expressions for $\ex_r(n,P_k)$ and $\ex_r(n,C_k)$ in Theorem \ref{main}.

(2) The proof of Theorem \ref{main} restricted to the case of $k$-paths is substantially simpler than the proof for $k$-cycles.

(3) It was recently shown by Bushaw and Kettle~\cite{BK} that the Tur\' an problem for 
disjoint $k$-paths can be easily solved once we know the extremal function for a single $k$-path.
As we have now solved the $k$-paths problem for all $r \ge 3$, the corresponding extremal questions for disjoint $k$-paths are 
also completely solved (for large $n$).  
A similar situation likely holds for disjoint $k$-cycles, as recently observed by Gu, Li and Shi~\cite{GLS}.

\subsection{The extremal function for minimal $k$-cycles and minimal $k$-paths.\mbox{ }} The related problems of determining ex$(n, \PP_k)$ and $\ex_r(n, \CC_k)$ have also received considerable attention, indeed the case of $\PP_2$ is the celebrated Erd\H os-Ko-Rado theorem. The last two authors~\cite{MV2} proved that $\ex(n,\PP_3)={n-1 \choose r-1}$ for all $r \ge 3$ and $n \ge 2r$.  The case of $\CC_3$ goes back to Chv\' atal~\cite{Chvatal} in 1973, and in \cite{MV} the last two authors proved that $\ex_r(n,\CC_3)={n-1 \choose r-1}$ for all $r \ge 3$ and $n \ge 3r/2$ thereby  settling an old conjecture of Erd\H os~\cite{E64}. They also proved some bounds for all $k,r$ and conjectured that both of these extremal functions are asymptotic to $\ell{n \choose r-1}$.
 F\"uredi, Jiang and Seiver~\cite{FJS} proved the conjecture in strong form and determined ex$(n, \PP_k)$ for all $k,r\ge 3$ and $n$ large.
 F\"uredi and Jiang~\cite{FJ} later determined  ex$(n, \CC_k)$ exactly for all $k \ge 3, r \ge 4$ and $n$ large.
 Our second theorem determines ex$_r(n, \CC_k)$ as well as the extremal  $\CC_k$-free $r$-graphs for all $r \geq 3$ and $n$ large.

\thm \label{loosemain}
Let $r \geq 3$, $k \geq 5$, and $\ell = \lfloor \frac{k-1}{2} \rfloor$. Then for sufficiently large $n$,
\[ \ex_r(n,\CC_k) = {n \choose r} - {n - \ell \choose r} +\left\{
\begin{array}{cl}
0 & \mbox{ if $k$ is odd,}\\
1 &\mbox{ if $k$ is even}
 \end{array}\right.
\]
with equality only for $r$-graphs of the form $S_{L}^r(n)$ with $|L|=\ell$ if $k$ is odd,
and $S_{L}^r(n)$ plus an edge when $k$ is even. Also for each $r \geq 3$,
\[ \ex_r(n,\CC_4) = {n \choose r} - {n - 1 \choose r} + \Bigl\lfloor\frac{n-1}{r}\Big\rfloor\]
with equality only for $r$-graphs of the form $S_L^r(n) \cup F$ where $F$ comprises $\lfloor \frac{n-1}{r}\rfloor$ disjoint edges.
\xthm

The proof is very similar to that of Theorem~\ref{main} and some steps are easier, so we only indicate the differences in the proofs.
The reader may observe that the approach also yields a proof for minimal paths that is substantially shorter than that in~\cite{FJS}.
Furthermore, we believe our methods with some additional refinements give polynomial bounds on $n$ relative to $r$ and $k$ above 
which Theorem \ref{main} and Theorem \ref{loosemain} hold. 

\subsection{The extremal problem for  Berge $k$-paths and $k$-cycles.\mbox{ }} Interesting results on the Tur\' an-type problems for Berge $k$-paths and Berge $k$-cycles,
were obtained by Bollob\' as and Gy\" ori~\cite{BG} and in a series of papers by Gy\" ori, Katona and Lemons, in particular, in~\cite{GL1,GL2,GKL}.
The bounds  differ from those in Theorems~\ref{main} and~\ref{loosemain}. In particular, they are linear in $n$ for $\ex_r(n,\mathcal{BP}_k)$.
We do not study $\ex_r(n,\mathcal{BC}_k)$ in this paper.
But if we forbid the family of Berge $k$-cycles or Berge $k$-paths in which no vertex belongs to at least $3$ edges,
then the answer is the same as in Theorem~\ref{loosemain}, apart from $k=4$:
the proof of the upper bound simply applies here, and the construction of $S_{L}^r(n)$ if $k$ is odd
 and $S_{L}^r(n)$ plus one edge if $k$ is even also applies. 
We remark that Tur\' an-type problems for Berge cycles with other additional restrictions have been extensively studied in the literature. 
Very recently, Jiang and Collier-Cartaino~\cite{JC} showed that a 2-linear $r$-graph on $n$ vertices with no $2k$-cycle has $O(n^{1 + 1/k})$ edges, 
generalizing the Even Cycle Theorem of Bondy and Simonovits~\cite{BS}.
 As another instance, for the minimal 4-cycle $C = \{e,f,g,h\}$ with $e \cup f = g \cup h$ and $e \cap f = g \cap h = \emptyset$,
 Erd\H{o}s~\cite{Erd2} conjectured $\ex_r(n,C) = O(n^{r - 1})$, and this was proved by F\"{u}redi~\cite{F} (see also~\cite{F,MV3,PV}).
 It seems likely that in this case the extremal $C$-free $r$-graphs for $r > 3$ are those in
  Theorem \ref{loosemain} for $k = 4$, and F\"{u}redi~\cite{F} conjectured $\ex_r(n,C) \sim {n - 1 \choose r - 1}$.

\subsection{Organization.\mbox{ }} We prove Theorem \ref{main} in four steps in Section \ref{proofofmain}; first we give
an asymptotic version, then a stability version followed by the proof of the exact result for cycles and the exact result for paths.
Theorem \ref{loosemain} is proved in Section \ref{proofofloosemain}. In Sections~\ref{lemmas}--\ref{randoms}
we prepare the background for passing from cycles and paths in the shadow of an $r$-graph to cycles and paths in the
$r$-graph itself.

\section{Notation and terminology}

\subsection{General notation.\mbox{ }}

Edges of an $r$-graph $H$ sometimes will be written as unordered lists, for instance, $xyz$ represents $\{x,y,z\}$.
For $X \subset V(H)$, let $H - X = \{e \in H : e \cap X = \emptyset\}$.
The {\em codegree} of a set $S = x_1 x_2 \dots x_s$ of vertices of $H$ is $d_H(S) = |\{e \in H : S \subset e\}|$; when $s=r-1$,
the {\em neighborhood in $H$ of $S$} is $N_H(S)=\{x: S \cup \{x\} \in H\}$, so that $|N_H(S)|=d_H(S)$.
For vertices $x,y$ in a hypergraph, an {\em $x,y$-path} is a path $P = e_0 e_1 \dots e_k$ where $x \in e_0 - e_1$ and $y \in e_k - e_{k - 1}$.

\subsection{Shadows in hypergraphs.\mbox{ }}

Now we state the crucial definitions involving shadows in hypergraphs.
Let $\partial H$ denote the $(r - 1)$-graph of sets contained in some edge of $H$ -- this is the {\em shadow} of $H$.
The edges of $\partial H$ will be called the {\em sub-edges} of $H$. If $G \subset \partial H$ and $F \subset H$ is obtained from
$G$ by adding distinct vertices of $V(H)- V(G)$ to each edge of $G$, then we say that $G$ {\em expands} to $F$.

For $2 \leq s < r$, let $\partial^1 H := \partial H$ and $\partial^s H = \partial^{s - 1} \partial H$.
The strategy to prove Theorem \ref{main} is to find a cycle in the shadow of an $r$-graph that can be
expanded to a cycle in the $r$-graph itself.

\defn Let $H$ be an $r$-graph. For $G \subset \partial H$ and $e \in G$, the {\em list} of $e$ is
\[ L_G(e) = N_H(e)- V(G).\]
The elements of $L_G(e)$ are called {\em colors}. We let $L_G = \bigcup_{e \in G} L_G(e)$ and
\[ \hat{G} = \{e \cup \{x\} : e \in G, x \in L_G(e)\}.\]
\xdefn

Note that all these definitions are  relative to the fixed host hypergraph $H$ and the fixed subgraph $G$ of $\partial H$.
A key idea is that if $C$ is a $k$-cycle or $k$-path in $\partial H$ and the family $\{L_C(e)\,:\,e \in C\}$ has a system of distinct representatives, then $\hat{C}$ contains a $k$-cycle or $k$-path, and so $H$ contains a $k$-cycle or $k$-path.

\section{Full, superfull and linear hypergraphs}\label{lemmas}

\subsection{Full subgraphs.\mbox{ }}
An $r$-graph $H$ is {\em $d$-full} if every sub-edge of $H$ has codegree at least $d$.
Thus $H$ is $d$-full exactly when the minimum non-zero codegree in $H$ is at least $d$.

The following lemma extends the well-known fact that any graph $G$ has a subgraph of minimum degree
at least $d + 1$ with at least $|G| - d|V(G)|$ edges.

\lem \label{fullsub}
For $r \geq 2,d \geq 1$, every $n$-vertex $r$-graph $H$ has a $(d + 1)$-full subgraph $F$ with
\[ |F| \geq |H| - d|\partial H|.\]
\xlem

\prf
A {\em $d$-sparse sequence} $S$ is a maximal sequence $e_1,e_2,\dots,e_m \in \partial H$ such that $d_H(e_1) \leq d$, and for all $i > 1$, $e_i$
 is contained in at most $d$ edges of $H$ which contain none of $e_1,e_2,\dots,e_{i - 1}$.
The $r$-graph $F$ obtained by deleting all edges of $H$ containing at least one member of a $d$-sparse sequence $S$ is $(d + 1)$-full. Since
 $S$ has length at most $|\partial H|$, we have $|F| \geq |H| - d|\partial H|$.
\xprf

\lem \label{full}
Let $r\geq 3$, $k \geq 3$ and let $H$ be a non-empty $rk$-full $r$-graph.  Then $C_k,P_{k-1} \subset H$.
\xlem

\prf
Consider the graph $F = \partial^{r-2} H$. Every edge of $H$ yields a $K_r$ in $F$, so
$F$ contains a $3$-cycle $C_3$. As $H$ is $rk$-full, each edge of $F$ is in at least $rk$ triangles in $F$.
We claim that $F$ contains a $k$-cycle: we start from $C_3$, and for $i=3,\ldots,k-1$,
obtain an $(i+1)$-cycle $C_{i+1}$ from $i$-cycle $C_i$ by using one of the at least $rk-i+2$ triangles containing an edge of $C_i$
and no other vertices of $C_i$. Let a $k$-cycle $C_k$ in $F$ have edges $f_1,\ldots, f_k$.
Choose in $H$ edges $e_1 = f_1 \cup g_1,\ldots ,e_k = f_k \cup g_k$ so that to maximize
the size of $Y = \bigcup_{i=1}^{k } e_i$. Suppose $C = \{e_1,\dots,e_{k}\}$ is not a $k$-cycle in $H$.
Then there are distinct $i,j$ such that $g_i \cap g_j \neq \emptyset$. Pick $v \in g_i \cap g_j$.
Let $Z = \{z \in V(H)\,:\, (f_i \cup g_i \cup \{z\}) - \{v\} \in H\}$. Since $H$ is $rk$-full, $|Z|\ge rk$.
As $C$ is not a $k$-cycle, $|Y| < rk$ and so there exists $z\in Z - Y$. Replacing $e_i$ with $e = (f_i \cup g_i \cup \{z\}) - \{v\}$,
we enlarge $Y$, a contradiction. So $H$ contains $C_k$ and thus $P_{k-1}$.
\xprf

\subsection{Superfull subgraphs.\mbox{ }}

\defn
An $\ell$-full $r$-graph $H$ is $\ell$-superfull if for every edge $e$ of $H$ at most one sub-edge of $e$ has codegree 
less or equal to $rk$.
\xdefn

\lem\label{doublinglemma}
Let $k,r \geq 3$, and let $H$ be an $\ell$-superfull $r$-graph such that $H$ contains a minimal $k$-cycle (respectively,
 a minimal $k$-path). Then $H$ contains a $k$-cycle (respectively, a $k$-path).
\xlem

\prf The proofs for paths and cycles are similar, so we only do the case of cycles. Let $C \subset H$ be a minimal $k$-cycle with  maximum $|V(C)|$.
If $C$ is not a $k$-cycle, then we find consecutive edges $f,g \in C$ with $|f \cap g| \geq 2$. Let $x,y \in f \cap g$.
Since $H$ is $\ell$-superfull, we may assume $d_H(f -\{x\})\ge rk$. Since $|V(C)| < rk$, we find $z \not \in V(C)$ such that $h = f \cup \{z\} - \{x\} \in H$.
 Then $C' = C \cup \{h\} - \{f\}$ has more vertices than $C$,
a contradiction.
\xprf

\smallskip

\lem\label{super}
Let $r \geq 3,k \geq 4$ and let $H$ be an $\ell$-superfull $r$-graph containing a set $W$ of at least $rk$ vertices
such that every $(r - 1)$-subset of $W$ has codegree exactly $\ell$. Let $G$ be the set of
all $(r - 1)$-subsets of $W$. If $H$ has no $k$-cycle or no $k$-path,
then for some set $L$ of $\ell$ vertices of $H - W$, $L_G(e) = L$ for every $(r - 1)$-set $e \subset W$.
\xlem

\prf If $e \cup \{x\} \in H$ for some $x \in W$, then all $(r-1)$-subsets of $e \cup \{x\}$ have codegree exactly $\ell$,
 contradicting the fact that $H$ is $\ell$-superfull.
Thus, $N_H(e) \cap W=\emptyset$  for all $e \in G$.

Suppose that $L_G(f)\neq L_G(e)$ for some  $e,f\in G$. Then there are  $e_1,e_2\in G$ such that
$|e_1\cap e_2|=1$ and $L_G(e_2)\neq L_G(e_1)$,
since from $|W| \geq rk \geq 4r$, for every two distinct $e,f\in G$, there is  $g\in G$ sharing exactly one vertex with each of $e$ and $f$.
In particular,
\begin{equation}\label{a1}
  |L_G(e_1)\cup L_G(e_2)|\geq \ell+1.
\end{equation}

{\bf Case 1: $\ell\geq 2$ and $H$ has no $k$-cycle.} Let $e_3,\ldots,e_{\ell+1}\in G$ be  such that $C = \{e_1,e_2,\ldots,e_{\ell+1}\}$ is an
$(\ell+1)$-cycle. By~\eqref{a1},
 the family $\{L_G(e_i) : 1 \leq i \leq \ell+1\}$ has a system of distinct
representatives $\{v_i \in L_G(e_i) : 1 \leq i \leq \ell+1\}$. As observed above,  $v_i \not\in W$ for all $i$.

Let $e_i \cap e_{i + 1} = \{w_{i+1}\}$ and $X_i = e_i \cup \{v_i\} - \{w_i,w_{i + 1}\}$, with subscripts modulo $\ell + 1$.
Then each of $X_i \cup \{w_i\}$ and $X_i \cup \{w_{i + 1}\}$
 has codegree at least $rk$ in $H$, since $H$ is $\ell$-superfull and $e_i$ has codegree exactly $\ell$.
Thus for each $1\leq i \leq \ell$, we can select
 edges $f_i,g_i \in H$ with $X_i \cup \{w_i\} \subset f_i$ and $X_i \subset \{w_{i + 1}\} \subset g_i$ forming a minimal $(2\ell + 2)$-cycle in $H$ if $k$ is even.
We let $f_{\ell+1} = g_{\ell+1} = e_{\ell+1}$ to
obtain a minimal $(2\ell + 1)$-cycle if $k$ is odd. In both cases, $H$ contains a minimal $k$-cycle, and so by Lemma~\ref{doublinglemma}, $H$ contains a $k$-cycle.

\medskip

{\bf Case 2: $\ell=1$  and $H$ has no $4$-cycle. } Let $e_3$ be a  sub-edge such that $\{e_1,e_2,e_3\}$ is a $3$-cycle.
For $i=1,2,3$, let  $L_G(e_i)=\{v_i\}$ and $e_i\cap e_{i+1}=\{w_i\}$. Note again  that $v_i \not\in W$. By symmetry, we may assume that $v_1\notin \{v_2,v_3\}$.
Since $H$ is $\ell$-superfull and $e_1$ has codegree exactly $\ell$, the sub-edges $e'=e_1-w_1+v_1$ and $e''=e_1-w_3+v_1$
 have codegrees at least $3r$. So we can select edges $g_1\supset e'$ and $g_2\supset e''$, so that
$\{e_2,e_3,g_1,g_2\}$ is a minimal $4$-cycle in $H$. Applying Lemma \ref{doublinglemma}, we conclude that $H$ contains a 4-cycle.

\medskip

{\bf Case 3: $H$ has no $k$-path.} We repeat Case 1, except we use an $(\ell+1)$-path instead of $C$.
\xprf

\subsection{Linear hypergraphs.\mbox{ }}

In the last two sections we showed how to pass from cycles and paths in the shadow of full and superfull subgraphs 
of an $r$-graph $H$ to
cycles and paths in $H$ itself. Here we consider the case that all sub-edges  have bounded codegrees.
The following fact is due to Erd\H os (see Theorem~1 in~\cite{E64}):

\prop[Erd\H os~\cite{E64}]\label{erd}
For $r,t\geq 2$ there exists $n_0=n_0(r,t)$ such that for all $n>n_0$, every $n$-vertex $r$-graph $H$ with  $|H|>n^{r-t^{1-r}}$
contains the complete $r$-partite $r$-graph $K^r_{t,\ldots,t}$.
\xprop

\defn
An $n$-vertex $r$-graph $H$ is {\em $(t,c)$-sparse} if every $t$-set of vertices lies in at most $c$ edges  of $H$.  If $c=1$, then  $H$ is $t$-linear.
\xdefn

The famous Ruzsa-Szemer\'edi $(6,3)$-Theorem~\cite{RS} shows that any linear $3$-graph on $n$ vertices and  $\Omega(n^2)$ edges contains $C_3$.
The following generalization was proved for $r = 3$ by S\'ark\H{o}zy and Selkow~\cite{SS}  using the Regularity Lemma.
We   avoid the use of regularity for $r > 3$:
\prop \label{linearbound}
Fix $c>0$ and $r,k \ge 3$. Let $H$ be an $n$-vertex $(r - 1,c)$-sparse $r$-graph not containing $P_k$ or not containing $C_k$.
Then $|H|=o(n^{r-1})$.
\xprop

\prf  It suffices to prove the result for $C_k$ since $P_k \subset C_{k+1}$.
In view of the S\'ark\H{o}zy--Selkow Theorem~\cite{SS}, we consider only $r\geq 4$.
Consider the graph with vertex set $H$ in which two vertices are adjacent if the intersection
of the corresponding edges of $H$
has  size $r-1$. Since $H$ is $(r-1,c)$-sparse, this graph has maximum degree less than $rc$,
so it contains an independent set $H_0$ of size at least $|H|/rc$. This means that $H_0$ is an $(r-1)$-linear $r$-graph.

Assume that $\epsilon>0$, $n$ is sufficiently large, and $|H_0|>\epsilon n^{r-1}$.
A standard averaging argument shows that there is an  $r$-partite subgraph of $H_0$ with at least  $(r!/r^r)|H_0|$ edges.
Let $X_1, \ldots, X_r$ be the $r$ parts and consider the edge-colored $(r-1)$-partite $(r-1)$-graph $H'\subset \partial H_0$
with parts $X_1, \ldots, X_{r-1}$ where the color of the edge $\{x_1, \ldots, x_{r-1}\}$, with $x_i \in X_i$ for $i \in [r-1]$ is the unique
$x_r \in X_r$ such that $\{x_1, \ldots, x_r\} \in H_0$. Such $x_r$ is unique as $H_0$ is $(r-1)$-linear.
 We will find a {\em rainbow} $C_k$ in $H'$ -- in other words a $k$-cycle in $H'$ whose lists have a system of distinct representatives. Since $|H'|>(\epsilon r!/r^r)n^{r-1}$ and $n$ is large, by
Proposition~\ref{erd}, there is a  complete $(r-1)$-partite $(r-1)$-graph
$K= K_{k,k,\ldots,k,s} \subset H'$ where $s=k^{2r-3}+1$ that has the same $(r-1)$-partition as $H'$.
Since $H_0$ is $(r-1)$-linear, every color class $S_c$ in $H'$ is   $(r-2)$-linear.
Now construct a hypergraph $H^*$ with vertex set $X_r$ (these are the colors of $H'$) and $s$ edges,
where the $i$th edge consists of the set of colors on edges incident to the $i$th vertex of $K$ in the part of size $s$.
Note that $H^*$ need not be uniform, but its edges have size at most $k^{r-2}$.

Pick a color $c$ (recall that $c$ is a vertex of $H^*$).  The number of edges of $H^*$ (these correspond to vertices of $K$ in $X_{r-1}$)
containing $c$ is at most $k^{r-2}$ since $S_c$ is $(r-2)$-linear. So $H^*$ has maximum degree at most $k^{r-2}$, edges of size at most  $k^{r-2}$, and size $s$.
Therefore $H^*$ has a matching $M$  of size  $s'=\lceil s/k^{2r-4}\rceil > k$ (by the greedy algorithm).  This means that
$K$ contains the complete $(r-1)$-partite $(r-1)$-graph $K'= K_{k,k,\ldots, k, s'}$ with partite sets $X'_1,\ldots,X'_{r-1}$, $|X'_1|=\ldots=|X'_{r-2}|=k$,
and $|X'_{r-1}|=s'$ (here $X'_{r-1}$ corresponds to $M$) such that
\begin{equation}\label{jul211}
\mbox{
  no two edges
$e,e'$ with the same color are incident to different vertices in  $X'_{r-1}$. }
\end{equation}

Let $x\in X'_{r-1}$. We claim that
\begin{equation}\label{jul21}
\mbox{there is a pair $\{e_1,e_2\}$ of edges in $K'$ of different colors such that $e_1\cap e_2=\{x\}$.}
\end{equation}
Indeed consider two edges $e=\{x_1,\ldots,x_{r-2},x\}$ and $e'=\{x_1,\ldots,x_{r-3},x'_{r-2},x\}$ of $K'$ that differ only in $(r-2)$th coordinate.
Since $H_0$ is an $(r-1)$-linear, they have different colors. Then for any edge $e''\in K'$ that shares only $x$ with $e\cup e'$,
either $\{e,e''\}$ or $\{e',e''\}$ satisfies~\eqref{jul21}.

Consider a $k$-cycle $C'=\{e_1,\ldots,e_k\}$ in $K'$ such that $e_1$ and $e_2$
satisfy~\eqref{jul21} and for every $i\neq 1$,  the vertex $v_i\in e_i\cap e_{i+1}$ is not in $X'_{r-1}$.
By~\eqref{jul211} and~\eqref{jul21}, $C'$ is a rainbow $k$-cycle in $K'$ and we expand it to a $k$-cycle in $H$. \xprf

\section{Cycles and paths from shadows}\label{coloringlemmas}

We now present the key lemmas which show how to expand $k$-paths and $k$-cycles in $\partial H$ to paths and
cycles in $H$ itself. Throughout this section, $r,k \geq 3$ and $\ell=\left\lfloor\frac{k-1}{2}\right\rfloor$.

\subsection{Paths.\mbox{ }}

\lem \label{easypaths}
Let $k \geq 3$, let $H$ be an $r$-graph and let $P = \{e_0,e_1,\dots,e_{2^{2\ell+1} - 1}\}$ be a $2^{2\ell+1}$-path in $\partial H$.
If $|L_P(e)| \geq \ell + 1$ for all $e \in P$, then $\hat{P}$ contains a $k$-path whose first edge contains $e_0$.
\xlem

\prf As $\lfloor (k - 1)/2 \rfloor = \lfloor (k - 2)/2\rfloor$ for $k$ even, it is enough to consider even $k \geq 4$. First we
prove the lemma for  $k = 4$, and then apply an inductive proof. The case $k = 4$ is split into two cases:

\medskip

{\bf Case 1: $L_P(e_0) \cap L_P(e_i) \neq \emptyset$ for some $i > 1$.}\\
Let $\alpha \in L_P(e_0) \cap L_P(e_i)$ and let $e_i,f,g,h \in P$ form a path vertex-disjoint from $e_0$ - this exists since $P$ has eight edges.
Define $L'(e) = L_P(e) - \{\alpha\}$ for $e\in P$. If we find distinct $\beta \in L'(f) $ and $\gamma \in L'(g) $,
then $\{e_0 \cup \{\alpha\},e_i \cup \{\alpha\},f\cup \{\beta\},g\cup \{\gamma\}\}$ is a 4-path.
Otherwise, $L_P(f) = L_P(g) = \{\alpha,\alpha'\}$ for some $\alpha'$. The same argument
with $f$ in place of $e_i$ shows $ L_P(g) = L_P(h) = \{\alpha,\alpha'\}$,
in which case the required 4-path is $\{e_0 \cup \{\alpha\},e_i \cup \{\alpha\},f \cup \{\alpha'\},h \cup \{\alpha'\}\}$.

\medskip

{\bf Case 2: $L_P(e_0) \cap L_P(e_i) = \emptyset$ for all $i > 1$.} \\
Let
$ L_P(e_0)=\{\alpha,\beta\}$. If $L_P(e_0) \cap L_P(e_1) \neq \emptyset$, say, $\beta\in L_P(e_1)$,
then by the case,
 we may pick distinct $\gamma \in L_P(e_2)$ and $\delta \in L_P(e_3)$
so that $\{e_0 \cup \{\alpha\},e_1 \cup \{\beta\},e_2\cup \{\gamma\},e_3 \cup \{\delta\}\}$
is a 4-path, as required.
Suppose $L_P(e_0) \cap L_P(e_1) = \emptyset$. If there is $\gamma\in L_P(e_1)\cap L_P(e_3)$,
then choose any $\lambda\in L_P(e_4)-\gamma$, and the edges
$e_0 \cup \{\alpha\},e_1 \cup \{\gamma\},e_3 \cup \{\gamma\},e_4 \cup \{\lambda\}$ form a 4-path.
Otherwise,
as $|L_P(e_i)| \geq 2$ for $i \geq 1$, we can choose all distinct $\alpha_1\in L_P(e_1)$, $\alpha_2\in L_P(e_2)$,
$\alpha_3\in L_P(e_3)$, and the  edges in the set $\{e_i\{\alpha_i\}\,:\,i=1,2,3\}$ together with
$e_0 \cup \{\alpha\}$ form a 4-path.

\medskip

Now suppose $k \geq 6$. If for some $i > 1$ we have $\beta \in L_P(e_0) \cap L_P(e_i)$,
let $P'=\{e_{i+1}, e_{i+2}, \ldots, e_{i+2^{k-3}}\}$ if $i \le 2^{k-3}+1$ and $P'=\{e_{i-1}, e_{i-2}, \ldots, e_{i-2^{k-3}}\}$ if $i> 2^{k-3}+1$ (note that $i-2^{k-3}\ge 2$). Let $e_0'=e_{i+1}$ if $i \le 2^{k-3} +1$ and $e_0'=e_{i-1}$ if $i>2^{k-3}+1$.
Let us remove $\beta$ from all lists of edges of $P'$.  Then $P'$ is a $2^{k-3}$-path all of whose lists have size at least $\ell$. So by induction on $k$,  $\hat{P} - \beta$ has a $(k - 2)$-path $\{f_2, f_3, \dots, f_{k - 1}\}$ where $e_0' \subset f_2$. Set
$f_0 = e_0 \cup \{\beta\}, f_1 = e_i \cup \{\beta\}$. Then $\{f_0, f_1, \dots,
f_k\}$ is the required $k$-path. So we may assume for all $i > 1$, $L_P(e_0) \cap L_P(e_i) = \emptyset$. If we find  $\gamma \in L_P(e_1)-L_P(e_0)$, then remove $\gamma$ from all lists $L_P(e_i)$ where $i \ge 2$. Let $\hat{P'} = \hat{P} - L_P(e_0) -\{\gamma\}$ if $\gamma$ exists and $\hat{P'}=\hat{P}-L_P(e_0)$ otherwise (in this case $L_P(e_1) \subset L_P(e_0)$).
By induction, $\hat{P'}$ contains a $(k - 2)$-path
$\{f_2,f_3, \dots, f_{k-1}\}$ with $e_2 \subset f_2$ as the lists sizes have reduced by at most one.
Set $f_0 = e_0 \cup \{\alpha\},f_1 = e_1 \cup \{\beta\}$ with $\alpha \neq \beta$, $\alpha \in L_P(e_0)$ and $\beta \in L_P(e_1)\cup \{\gamma\}$ (if $\gamma$ exists  we may choose $\beta=\gamma$); this works since $|L_P(e)| \geq 2$ for $e \in P$.
Now $\{f_0, f_1, \dots, f_{k-1}\} \subset \hat{P}$ is a $k$-path. \xprf

\subsection{Cycles.\mbox{ }}

To extend Lemma \ref{easypaths} to $k$-cycles, we need the following technical definition.

\defn\label{c1}
Let $H$ be an $r$-graph where $r \geq 3$. Let $\Psi_t(H)$ be the set of complete $(r - 1)$-partite $(r - 1)$-graphs
$G \subset \partial H$ with parts of size $t$ and $|L_G(e)| > \ell$ for all $e \in G$, and if $r = 3$ and $k$ is odd, then in addition
for $xy \in G$, there is $xy\alpha \in \hat{G}$ such that

\quad $(a)$   $\min\{d_H(x\alpha),d_H(y\alpha)\}\geq 2$  and

\quad $(b)$ $\max\{d_H(x\alpha),d_H(y\alpha)\}\geq 3k+1$.
\xdefn

The additional technical conditions for $r = 3$ and $k$ odd will become apparent in the proof of Case 2 of Lemma \ref{cyclelemma} below.
We also will use the following   consequence of Hall's Theorem:

\lem \label{hall}
Let $p \geq 1$ and $q \in \{2p,2p + 1\}$, and let $S_1,S_2,\dots,S_q$ be sets such that $S_i \cap S_j = \emptyset$
for $i \leq p$ and $j \geq  p + 2$, and $|S_i| > p$ for $i \leq p$ and $|S_i| \geq p$ for $i > p$. Then
$\{S_1,S_2,\dots,S_q\}$
has a system of distinct representatives, unless $q = 2p + 1$ and all $S_j$ for $j > p$ are all equal and of size $p$.
\xlem

\prf If the lemma is false, then by Hall's Theorem, there is
$I \subset [q]$ such that $|\bigcup_{i \in I} S_i| < |I|$. As $S_i \cap S_j = \emptyset$
for $i \leq p$ and $j \geq p + 2$,  $I \subset [p + 1]$ or $I \subset [p + 1,q]$.
It is not possible that $I \subset [p + 1]$, since $|S_i| > p$ for $i \leq p$. If $I \subset [p + 1,q]$,
then since $|S_i| \geq p$ for $i \in I$, the only possibility is $q = 2p + 1$ and $I = [p + 1,q]$ and $|\bigcup_{i \in I} S_i| = p$.
In this case all $S_i$ for  $i \in I$ are identical.
\xprf

\lem\label{cyclelemma} Let $r \geq 3,k \geq 4$, and let $H$ be a $C_k$-free $r$-graph. If $t$ is large enough then $\Psi_t(H) = \emptyset$.
\xlem

\prf Suppose $G \in \Psi_t(H)$.
Let $M$ be a set of $s = 2^{k - 2}(r - 1)$ pairwise disjoint edges of $G$. If there exists $\alpha \in L_G(e)$ for all $e \in M$,
let $F \subset G$ be a complete $(r - 1)$-partite subgraph of $G$ with $V(F) \subset V(M)$, $|f \cap V(M)| = 1$ for all $f \in F$, and parts of size $2^{k-2}$. We show that $\hat{F}$ contains a $(k - 2)$-path avoiding $\alpha$.
For $k \ge 5$, $F$ contains a $2^{k-2}$-path, so by Lemma \ref{easypaths}, $\hat{F}$ contains a $(k-2)$-path.
If $k=4$ and $F$ has  lists of size $1$ after removing $\alpha$, we cannot use Lemma \ref{easypaths} to find a $(k-2)$-path as $k-2 < 3$.
To find a 2-path in $\hat{F}$ in this case, consider any $3$-path $\{f_1,f_2,f_3\}$ in $F$.
 Suppose $\beta_i\in L_G(f_i)-\alpha$ for $i=1,2,3$. If $\beta_1=\beta_3$, then $\{f_1\cup {\beta_1}, f_3\cup {\beta_1}\}$
 is a $2$-path; otherwise either $\{f_1\cup {\beta_1}, f_2\cup {\beta_2}\}$ or
 $\{f_2\cup {\beta_2}, f_3\cup {\beta_3}\}$ is a $2$-path. For all $k \ge 4$  we have found  $x,y \in V(F) \subset V(M)$
and an $xy$-path $\hat{P} \subset \hat{F} - \{\alpha\}$ of length $k - 2$.
 Picking edges $e,f \in M$ with $x \in e$ and $y \in f$, $\hat{P} \cup \{e \cup \{\alpha\},f \cup \{\alpha\}\}$ is a $k$-cycle in $\hat{G}$, a contradiction.
We conclude that
\begin{equation}\label{nocolor}
\mbox{no color appears in the lists of }s\mbox{ pairwise disjoint edges of }G.
\end{equation}
For every $e\in G$, fix  a subset $L'_G(e)$  of $L_G(e)$ with $|L'_G(e)|=\ell+1$.
Let $m = \lfloor t/(s + 2)\rfloor$.
For $i \in [m]$, let $F_i \subset G$ be vertex-disjoint complete $(r - 1)$-partite graphs with parts of size $s + 2$, and $L'_i = \bigcup\{L'_G(e) : e \in F_i\}$.
Then $|L'_1| \leq (\ell + 1)|F_1| < (s + 2)^r$.  For each color $\alpha \in L'_1$, by (\ref{nocolor}),
there are at most $s$ different $i$ for which  $\alpha \in L'_i \cap L'_1$. So $L'_i \cap L'_1 \ne\emptyset$ for at most $(s + 2)^{r+1}$ values $i \in [m]$.
Choose $t$ so that $m > (s + 2)^{r+1}$. Then for some $i > 1$, $L'_i \cap L'_1 = \emptyset$, say for $i = 2$. Let $F = F_1 \cup F_2$ and let $X,Y$ be two parts of $F$.
Select $e \in G$ with $e \cap V(F_1) = \{x\} \subset X$ and $e \cap V(F_2) = \{y\} \subset Y$.

\medskip

{\bf Case 1: $r > 3$, or $r = 3$ and $k$ is even.} Let $e \cup \{\alpha\} \in \hat{G}$.
By the symmetry between $L'_1$ and $L'_2$ we may suppose $\alpha\not \in L'_1$.
Let $q = k - 1, p = \ell$ and let $U$ be a part of $F_1 - \{x\}$  and $V$ be a different part in $F_2 - \{y\}$.
Let $f$ be
 any edge $f \in G$ with $|f \cap U| = 1 = |f \cap V|$ and $|f \cap V(F)| = 2$.
Since $U$ and $V$ are subsets of different parts in $F$ and $r > 3$, or $r = 3$ and $k$ is even, there is a
$q$-path $Q = \{f_1,f_2,\dots,f_q\}$ from
$x$ to $y$ in $G$ with $f_i \subset F_1$ for $i \leq p$, $ f_{p + 1}=f$, and $f_i \subset F_2$ for $i > p + 1$.
 If $Q$ expands to a $q$-path $\hat{Q} \subset \hat{G} - \alpha$, then $\hat{Q} \cup \{e \cup \{\alpha\}\}$ is a
$k$-cycle in $\hat{G}$, a contradiction. Therefore
\begin{equation}\label{noexpand1}
Q\mbox{ does not expand to a }q\mbox{-path in }\hat{G} - \alpha.
\end{equation}
Now let $S_i = L'_G(f_i) - \alpha$ for $1 \leq i \leq q$. Since $L'_1 \cap L'_2 = \emptyset$, we have $S_i \cap S_j = \emptyset$ for $i \leq p$ and $j > p + 1$,
and since $\alpha \not \in L'_1$, $|S_i| > p$ for $i \leq p$, and $|S_i| \geq |L'_G(f_i)| - 1 \geq p$ for $i > p$.
By (\ref{noexpand1}), the family $\{S_1,S_2,\dots,S_q\}$ has no system of distinct representatives.
By Lemma \ref{hall}, all $S_i$ for $i > p$ are identical of size $p=\ell$, and since $|L_g'(f_i)|=\ell+1$, we have  $\alpha \in L'_G(f)$. Since $f$ was any edge with
$|f \cap U| = 1 = |f \cap V|$ and $|f \cap V(F)| = 2$, $G$ is complete $(r-1)$-partite, $t$ is large, and $|U|,|V| \geq s$,
 we have $s$ disjoint edges of $G$ whose lists all
 contain $\alpha$, contradicting (\ref{nocolor}). This finishes Case 1.

\medskip

{\bf Case 2: $r = 3$ and $k$ is odd.} Let $q = k - 2$ and $p = \ell - 1$, so $q = 2p + 1$. Since $G\in \Psi_t(H)$, some
 $xy\alpha \in \hat{G}$ satisfies (a) and (b) in  Definition~\ref{c1}. Again, since
 $L'_1\cap L'_2=\emptyset$, we may suppose $\alpha\not \in L'_1$.
 By symmetry we may assume $d_H(x\alpha) > 3k$ and $d_H(y\alpha) > 1$. Choose an edge $y\alpha\beta \in H$ with $\beta \neq x$.
Note that possibly $\beta \in V(G)$. For $i=1,2$, let $X_i=X\cap V(F_i) - \{x,\beta\}$ and $Y_i=Y\cap V(F_i) - \{y,\beta\}$.
Let
 $f \in G$  be such that
\[
|f \cap X_1| = 1 = |f \cap Y_2|  \mbox{ if }  q\equiv 1 \mbox{ (mod 4),} \quad |f \cap X_2| = 1 = |f \cap Y_1|  \mbox{ if } q\equiv 3 \mbox{ (mod 4).}
\]
  Since $q$ is odd, there is a
$q$-path $Q = \{f_1,f_2,\dots,f_q\}$ from $x$ to $y$ in $G$ with $f_i \subset F_1$ for $i \leq p$, $f_{p + 1}=f$, and $f_i \subset F_2$ for $i > p + 1$.
If $Q$ expands to a $q$-path $\hat{Q} \subset \hat{G} - \alpha - \beta$, then select $\gamma \in V(H) - V(\hat{Q} )-\alpha-\beta$
so that $x\alpha\gamma \in H$ -- this is
possible since $d_H(x\alpha) > 3k$ -- and then $\hat{Q} \cup \{x\alpha\gamma,y\alpha\beta\}$ is a
$k$-cycle in $\hat{G}$. So
\begin{equation}\label{noexpand2}
Q\mbox{ does not expand to a }q\mbox{-path in }\hat{G} - \alpha - \beta.
\end{equation}
Let $S_i = L'_G(f_i) - \alpha - \beta$.
Since $L'_1 \cap L'_2 = \emptyset$, we have $S_i \cap S_j = \emptyset$ for $i \leq p$ and $j > p + 1$,  and since $\alpha \not \in L'_1$,
$|S_i| = |L'_G(f_i) - \beta| \geq \ell > p$ for $i \leq p$, and $|S_i| \geq |L'_G(f_i)| - 2 \geq p$ for $i > p$. By (\ref{noexpand2}), the family
 $\{S_1,S_2,\dots,S_q\}$ has no system of distinct representatives.
By Lemma~\ref{hall}, all $S_i$ for $ i > p$ are identical, and in particular, $\alpha \in L'_G(f)$.
Since $f$ was an arbitrary edge joining $X_1$ to $Y_2$ or joining $X_2$ to $Y_1$ and $|X_i|, |Y_i| \geq s$ for $i=1,2$, this contradicts (\ref{nocolor}).
\xprf

\section{Random sampling}\label{randoms}

We use a random sampling technique and Lemmas~\ref{cyclelemma} and~\ref{easypaths}
to find $k$-cycles and $k$-paths
in an $r$-graph $H$  when $H$
 has many sub-edges of codegree at least $\ell + 1$.

\lem \label{kttlemma}
Let $\delta >0$, $r\geq 3$ and $k\geq 4$. Let $H$ be an  $r$-graph, and
$E \subset \partial H$ with $|E| > \delta n^{r - 1}$.  Suppose that $d_H(f) \ge \ell+1$  for every $f \in E$
and, if $r=3$ and $k$ is odd, then in addition, for every $f=xy\in E$
there is $e_f=xy\alpha\in H$ such that $\min\{d_H(x\alpha),d_H(y\alpha)\}\geq 2$  and
 $\max\{d_H(x\alpha),d_H(y\alpha)\}\geq 3k+1$. Then for large enough $n$, $H$ contains $P_k$ and $C_k$.
 \xlem

\prf By Lemmas~\ref{cyclelemma}
and~\ref{easypaths}, it is enough to prove that $\Psi_t(H)\neq\emptyset$ for a large enough $t$.

Let $m=\ell+1$ and $T$ be a random subset of $V(H)$ obtained by picking each vertex independently with probability $p=1/2$. Let
\[ F = \{f \in E : f \subset T, |N_H(f)- T| \geq m, e_f-f \not\subset T\}.\]
For $f \in  E$ and any choice of edges $e_1,e_2,\dots,e_{m} \in H$ containing $f$ such that $e_1 =e_f$, the probability that
$f \subset T$ and $e_i-f \not\subset T$ for $i \in [m]$ is exactly $p^{r - 1}(1 - p)^m$. Therefore
$$\mathbb E(|F|) \ge |E| p^{r-1}(1-p)^{m}\ge \delta 2^{-m-r+1}n^{r-1}.$$
So there is a $T \subset V(H)$ with $|F| \geq \delta 2^{-m-r+1}n^{r - 1}$.
If $n$ is large enough, Proposition~\ref{erd} gives a complete $(r - 1)$-partite $G \subset F$
with parts of size $t$. Since $|L_G(f)| \geq |N_H(f)- T| \geq m$ for $f \in G$, $G\in \Psi_t(H)$ for $r\geq 4$ and for even $k$
when $r=3$. Suppose $r=3$ and $k$ is odd. Then since for every $f\in G$, $e_f\in \hat{G}$, again $G\in \Psi_t(H)$.
\xprf

\section{Proof of Theorem \ref{main}}\label{proofofmain}

\subsection{Part I : Asymptotics.\mbox{ }}\label{asym}

\thm \label{asymptotics}
 Let $r \geq 3,k \geq 4$.\\
 $(a)$  If $H$ is an $n$-vertex $(\ell + 1)$-full $r$-graph and $C_k \not \subset H$ or $P_k \not \subset H$, then $|H| = o(n^{r - 1})$. \\
 $(b)$  $\ex_r(n,P_k) \sim\ex_r(n,C_k) \sim \ell {n \choose r - 1}$.
\xthm

\prf To prove (a), we first show
\begin{equation}\label{jul1}
|\partial H| = o(n^{r-1}).
\end{equation}
Suppose  that $|\partial H| > \delta n^{r - 1}$ where $\delta > 0$, and $n$ is large.
If $r > 3$ or $r=3$ and $k$ is even, then by Lemma \ref{kttlemma} with $E=\partial H$, if $t$ is large enough, then
$H$ contains a $k$-cycle and  a $k$-path,
a contradiction.

 For $r = 3$ and $k$ odd, let $H^*$ be the set of edges of $H$
containing no pair of codegree at least $3k$. Then $H^*$ is $(2, 3k)$-sparse, so by Proposition \ref{linearbound}, $|H^*| = o(n^{2})$. Let $F = \partial H -\partial H^*$ so that for every $f \in F$,
there is an edge $e \in H$ containing $f$ and containing a pair $f'$ with $d_H(f') > 3k$  (possibly, $f'=f$).
 Then $|F| \geq |\partial H| - |\partial H^*| \geq \delta n^{2} - o(n^{2})> (\delta/2)n^{2}$ if $n$ is large enough.

If all edges of $H$ containing a pair $f \in F$ have all their sub-edges of codegree greater than $3k$,
map $f$ to itself. Otherwise, pick an edge of $H$ containing $f$ and containing some pair $f'$ of codegree at most $3k$, and
map $f$ to $f'$ (again $f = f'$ is possible).
This map is at most $6k$ to one, and therefore we have a set $E$ of $(\delta/12k)n^2$ pairs in $\partial H$ each
of codegree at least $\ell + 1$ in $H$ and each $f \in E$ is contained in some edge  $e_f \in H$ in which some other pair has
codegree at least $3k + 1$. Since $H$ is $(\ell+1)$-full, the conditions of Lemma \ref{kttlemma} hold for $E$, and
so $H$ contains a $k$-cycle and  a $k$-path, a contradiction.
So we proved~\eqref{jul1} in both cases.

Now by Lemma \ref{fullsub}, $H$ has an $r(k+1)$-full subgraph $H'$ with
\[ |H'| \geq |H| - r(k+1)|\partial H|.\]
By Lemma \ref{full}, if $H' \neq \emptyset$, then $P_k,C_k \subset H' \subset H$, which is a contradiction. we conclude $H' = \emptyset$, and
so $|H| \leq r(k+1)|\partial H| = o(n^{r - 1})$, which proves (a).

Now we determine the asymptotic value of $\ex_r(n,C_k)$ and $\ex_r(n,P_k)$.
The construction $S_{L}^r(n)$ in the statement of Theorem~\ref{main} shows
$\ex_r(n,C_k),\ex_r(n,P_k)\geq {n \choose r} - {n - \ell \choose r}\sim \ell {n \choose r - 1}$.
Suppose  $H$ is an $r$-graph and $C_k \not \subset H$ or $P_k \not \subset H$.
By Lemma \ref{fullsub}, $H$ has an $(\ell + 1)$-full subgraph $H'$ with
$|H'|\geq |H| - \ell|\partial H|$. By (a),
$|H'|=o(n^{r - 1})$. So $|H|\leq |H'| + \ell|\partial H| \leq o(n^{r - 1})+\ell {n \choose r - 1}$.
\xprf

\subsection{Part II : Stability.\mbox{ }}\label{stability}

\thm \label{stab}
Fix $r \geq 3$, $k\geq 4$ and let $H$ be an $n$-vertex $r$-graph  with $|H| \sim \ell {n \choose r - 1}$
containing no $k$-cycle or no $k$-path.
Then there exists $G^* \subset \partial H$ with $|G^*| \sim{n \choose r - 1}$ and a set $L$ of $\ell$ vertices of $H$ such that $L_{G^*}(e) = L$ for every $e\in G^*$.  In particular, $|H-L|=o(n^{r-1})$.
\xthm

\prf Let $H^*$ be the set of edges of $H$ not containing any sub-edge of codegree at least $rk + 1$. Then $H^*$ is $(r-1, rk)$-sparse, so
Proposition~\ref{linearbound} implies  $|H^*|=o(n^{r-1})$. Let $H' = H  - H^*$, so $|H'| \sim |H|$.
We construct sequences $f_1,f_2,\ldots,f_q \in \partial H'$ and
 $H_0,H_1,\ldots,H_q \subset H$ with $H_0 = H'$ as follows.
Suppose $H_i$ is constructed and let $d_i(f)=d_{H_i}(f)$. A sub-edge $f$ of $H_i$ is of type
\begin{center}
\begin{tabular}{cp{5.8in}}
(i) & if $d_i(f) < \ell$, \\
(ii) & if $d_i(f) = \ell$ and some $e \in H_i$ containing $f$ contains a sub-edge $g \neq f$ with $d_i(g) = \ell$, \\
(iii) & if $\ell < d_i(f) < rk$.
\end{tabular}
\end{center}
If $H_i$ has no sub-edges of types (i) -- (iii), let $q = i$ and stop.
Otherwise, let $f$ be a sub-edge of $H_i$ of minimum type, and $H_{i+1} = H_i - \{e \in H_i : f \subset e\}$ and $f_{i+1} = f$.

Every sub-edge $f \in \partial H_q$ has $d_q(f)\ge \ell$ (since $f$ is not type (i)) so $H_q$ is certainly $\ell$-full. Also, no edge has more than one sub-edge of codegree less than $rk$, for then we have a sub-edge of type (ii) or (iii). Therefore $H_q$ is $\ell$-superfull.

\smallskip

{\bf Claim 1.} $|\partial H_q| \sim{n \choose r-1}$. \\
{\it Proof.} Let $E$ be the set of $f_i$ of type (iii), and for each $f \in E$, let $e_f$ be any edge of $H'$ containing $f$.
Suppose $|E| > \delta n^{r - 1}$. If $r\geq 4$ or $r=3$ and $k$ is even, this contradicts
 Lemma~\ref{kttlemma}. Let $r=3$ and $k$ be odd. By definition every edge of $H_i$
 containing  $f_i$
 of type (iii) has each of its subedges of codegree at least $\ell \geq 2$ and $d_H(f_i) \geq \ell + 1$.
 Since every edge in $H'$ contains some pair of
 codegree at least $3k + 1$ in $H$, the conditions of Lemma~\ref{kttlemma} are met by $E$. Again, by this lemma, $H$ contains $P_k$ and $C_k$, a contradiction.
 So, $|E| = o(n^{r - 1})$. Since we have deleted $q$ sub-edges, $|\partial H_q| \leq {n \choose r - 1} - q$.
Note that if a sub-edge of type (ii) was chosen, then $H_{i+1}$ will have a sub-edge of type (i).
So, if $\epsilon>0$  and  $q =\epsilon \binom{n}{r - 1}$, then for $n$ sufficiently large,
 $$|H_q| \geq |H'| - q(\ell-\frac{1}{2})- rk|E| \geq \ell |\partial H_q| -o(n^{r-1})+ \frac{\epsilon}{2}\binom{n}{r - 1} - rk|E| \geq \ell |\partial H_q| + \frac{\epsilon}{4}\binom{n}{r - 1}.$$
 By Lemma \ref{fullsub}, $H_q$ has an $(\ell + 1)$-full subgraph with at least $\frac{\epsilon}{4}\binom{n}{r - 1}$ edges,
 contradicting Theorem \ref{asymptotics}.  So $q=o\left(n^{r-1}\right)$, and $\ell |\partial H_q|\le |H_q|\le \ell|\partial H_q|+o(n^{r-1})$, which imply $|\partial H_q| \sim {n \choose r-1}$. \qed

Let $G'$ be the subgraph of $\partial H_q$ formed by the sub-edges
of codegree $\ell$ in $H_q$.
\medskip

{\bf Claim 2.} {\em $|G'| \sim{n \choose r-1}$.} \\
{\it Proof.} Let $G''=\partial H_q-G'$. Since $H_q$ is $\ell$-superfull, the codegree
of every $f\in G''$ is at least $\ell+1$. So if $r\geq 4$ or $r=3$ and $k$ is even, then by Lemma~\ref{kttlemma} with
$E=G''$, $|G''|=o(n^{r-1})$. If $r=3$ and $k$ is odd, then $\ell\geq 2$ and since $H_q$ is $\ell$-superfull,
the conditions of Lemma~\ref{kttlemma} are satisfied. So again we get
 $|G''|=o(n^{r-1})$, and thus
 $|G'| \sim{n \choose r-1}$ as required.  \qed

{\bf Claim 3.} {\em For each $rk$-clique $K \subset G'$, there exists $L \subset V(H_q) \backslash V(K)$ with $|L| = \ell$ and $L_K = L$.}  \\
{\it Proof.} As $H_q$ is $\ell$-superfull, this follows from Lemma \ref{super}. \qed

{\bf Claim 4.} {\it For some $G^* \subset G'$, $|G^*| \sim {n \choose r - 1}$ and all edges
of $G^*$ have the same list in $H_q$.}\\
{\it Proof.} Let $N$ be the number of $rk$-cliques in $G'$. Since $|G'| \sim {n \choose r - 1}$,
we easily see that $N \sim {n \choose rk}$.
By averaging, some edge $e^* \in G'$ is contained in at least
\[ \frac{N}{|G'|}{rk \choose r - 1}\]
$rk$-cliques in $G'$.

Since $H_q$ is superfull, if $K$ is an $rk$-clique in $G'$ containing
$e^*$ and $\alpha \in L_K(e^*)$, then for every $v\in e^*$, the sub-edge
$e^* + \alpha - v$ has codegree more than $rk > \ell$, and hence is not in $G'$. Thus
$L_K(e^*) \cap V(K') = \emptyset$ for every two $rk$-cliques $K,K' \subset
G'$ containing $e^*$. We stress that the lists here are taken
in $H_q$. In particular, there exists a set of $\ell$ vertices $L \subset V(H_q)$ such that $L_K(e^*) = L$
for every $rk$-clique $K \subset G'$ containing $e^*$. Let $G^* \subset G'$ be the set of edges of $G'$ contained in a common
$rk$-clique of $G'$ with $e^*$. By Claim 3, $L_K(f) = L$ for all $f \in G^*$. The number of pairs
$(K,f)$ where $K$ is an $rk$-clique in $G'$ containing $e^*$ and $f \in K$ is disjoint from $e^*$ is at least
\[ \frac{N{rk \choose r - 1}{rk - r + 1 \choose r - 1}}{|G'|}.\]
The number of $rk$-cliques containing both $e^*$ and $f$ is at most ${n \choose rk - 2r + 2}$. We conclude
\[ |G^*| \geq \frac{N{rk \choose r - 1}{rk - r + 1 \choose r - 1}}{|G'|{n \choose rk - 2r + 2}}.\]
Using $|G'| \sim {n \choose r - 1}$ and $N \sim {n \choose rk}$, a straightforward calculation shows
$|G^*| \sim {n \choose r- 1}$. \xprf

\subsection{Part IIIa : Exact result for cycles.\mbox{ }}\label{exact}

Fix $r \ge 3$, $k \ge 4$ and let  $n$ be large.
Let $H$ be an $n$-vertex $r$-graph containing
no $k$-cycle and with $|H| = {n \choose r} - {n - \ell \choose r} + f(n,k,r)$,
where  $f(n,k,r)=0$ if $k$ is odd, $f(n,k,r)=\mbox{ex}_r(n-\ell, \{P_2, 2P_1\})={n-\ell -2 \choose r-2}$ if $k$ is even and $(k,r)\ne (4,3)$ and $f(n,4,3)=\mbox{ex}_3(n-\ell, P_2)$.

Let $\beta=1/10$. Theorem~\ref{stab} implies that
for $n$ sufficiently large, $\ex_r(n, C_k) < 2 \ell {n \choose r-1}$ and consequently,
there is a  $c=c(k,r)$ such that $\ex_r(n, C_k)< cn^{r-1}$ for all $n \ge 1$.  Choose $\alpha$ sufficiently small so that
\begin{equation} \label{alpha}
c 2^{r-1} (k^3 r^r)^{r-1}\alpha^{(r-2)}< \beta/2.
\end{equation}
Finally, choose $n$ sufficiently large so that all inequalities involving $\alpha, k,r$ in the proof below are valid. By Theorem~\ref{stab}, there exists $L = \{x_1, \dots,x_{\ell}\} \subset [n]$ such that
$|H - L| \leq \alpha n^{r-1}$. Let $B=H-L$ be the set of edges of $H$ that are disjoint from $L$ so $|B|<\alpha n^{r-1}$. If $k$ is odd, then we shall show that
 $B = \emptyset$. If $k$ is even  then we shall show that $B$ is an extremal family with no $P_2$ and $2P_1$ unless $k=4, r=3$,
 in which case $B$ is an extremal family with no $P_2$. This proves both the extremal result and the characterization of equality. Let
$$M = \left\{e\in {[n] \choose r}- H: e \cap L \neq \emptyset\right\},$$
so that $$|B| = |M|+ f(n,k,r).$$
  If $M = \emptyset$, then we are done, so we may suppose for a contradiction that $M \neq \emptyset$ and $|B|>f(n,k,r)$. Set $m:=|M|$ so that
$m \le |B|< \alpha n^{r-1}$. 

\medskip

{\bf Claim 1.} {\it There exist pairwise disjoint $(r - 2)$-sets $Z_1,Z_2,\dots,Z_{kr} \subset V(H)- L$
such that for each $i \in [kr]$ and $j \in [\ell]$
$$\hbox{d}_{H}(Z_i \cup \{x_j\})\ge
n-r+1-\frac{krm}{{n -\ell \choose r-2}}.$$
If $r\ge 4$  there exists an additional $(r-2)$-set
$Z_{kr+1}$ that is disjoint from $Z_i$ for $i \in [kr-1]$
and $|Z_{kr+1} \cap Z_{kr}|=1$
}

{\it Proof.} Pick an $(r-2)$-set $T \subset V(H)- L$ uniformly at random. Let $\overline{H} = \{e \subset V(H) : |e| = r, e \not \in H\}$.
For $j \in [\ell]$, let $$X_j=d_{\overline{H}}(T \cup \{x_j\})=n-r+1-d_{H}(T \cup \{x_j\}).$$
 In other words, $X_j$ counts the number of $r$-sets $e \not\in H$ with $T \cup \{x_j\} \subset e$.
The number of $r$-sets $e \supset \{x_j\}$ with  $e \not\in H$ is at most $m$.
 For each such $e$, let $X_j(e)$ be the indicator for the event that $T \subset e$. Then
$$\mathbb E(X_j)=\sum_{e} \mathbb E(X_j(e))\le m \frac{{r-1 \choose r-2}}{{n-\ell \choose r-2}} < \frac{rm}{{n-\ell \choose r-2}}.$$
By Markov's inequality,
$$\mathbb P\left(X_j>\frac{k rm}{{n -\ell \choose r-2}}\right)<1/k.$$   This implies that
$$\mathbb P\left(\exists j: X_j> \frac{krm}{{n -\ell\choose r-2}}\right)< \ell/k<1/2.$$ In other words, the number of $T$ for which $d_{H}(T \cup \{x_j\})
\ge n-r+1-krm/{n-\ell \choose r-2}$ for all $j$ is at least ${n-\ell \choose r-2}/2$.

Now consider the family of all $(r-2)$-sets
 described above, and let $T_1, \ldots, T_t$ be a maximum matching in this family.
If $t<kr$, then all other sets of this family have an element
within $\cup_iT_i$, which implies that the number of such $T$ is less than
${n -\ell\choose r-2}/2$, because $n$ is sufficiently large. This  contradiction shows that $t \ge kr$.

If $r\ge 5$, then by a result of Frankl~\cite{Frankl} that
 $\mbox{ex}_{r-2}(n-\ell, P_2)=
O(n^{r-4})$, we can find two sets $T_1$, $T_2$
 with $|T_1 \cap T_2|=1$ and then find the remaining $kr-1$ sets using  the greedy procedure described above.
If $r=4$, then we use the fact that a graph with $\Omega(n^2)$ edges has a 2-path together with a disjoint from it matching of size $kr-1$. \qed

\medskip

{\bf Claim 2.} {\it Let  $Z=\cup_i Z_i$ and $Y = V(H) -(L \cup Z)$. Then there exists a set
$D \subset Y$ such that $H$ contains all edges of the form $Z_i \cup \{x_j,y\}$, for all $i \in [kr]$, $x_j \in L$ and $y \in D$ and
\[ |D| = n - \ell k r - \Big\lceil\frac{k^3r^r m}{n^{r-2}}\Big\rceil.\]}

{\it Proof.} For each $i\in [kr]$ and $j \in [\ell]$,  let
$S_{i,j}=\{y \in Y: Z_i \cup \{x_j,y\} \not\in H\}$.
Claim 1 implies that $|S_{i,j}|<krm/{n-\ell \choose r-2}$. Let $S=\cup_{i,j}S_{i,j}$. Then
$$|S|< \frac{(kr \ell) krm}{{n -\ell\choose r-2}}<  \frac{k^3r^r m}{n^{r-2}}.$$
We may add points arbitrarily to $S$ till $D:=Y-S$ has the required size.
   \qed

\medskip

{\bf Claim 3.} {\it No two edges $e,e' \in B$ have $|e \cap e'| = 1$ and $(e - e') \cap D \neq \emptyset$ and
$(e' - e) \cap D \neq \emptyset$.
If $k\ge 5$ is odd, then no edge $e \in B$ has $|e \cap D| \geq 2$. If $k\ge 6$ is even and $r=3$, then there are no two disjoint edges each with at least two points in $D$.
}

{\it Proof.} For $k$ even and $|e \cap e'|=1$ suppose $u \in e - e'$ and $v \in e' -e$. Then there is a path $P$ of
length $k - 2$ in $H$ between $u$ and $v$  consisting of edges $Z_i \cup \{x_j,y\}$ with $y \in D$
and such that $V(P) \cap (e \cup e') = \{u,v\}$. All vertices of $L$ will have degree two in $P$. Now $P \cup \{e,e'\}$ is a $k$-cycle in $H$.
 For $k\ge 5$ odd and $r \ge 4$, we repeat the same argument except that we use $Z_{kr-1}$ and $Z_{kr}$ which have a common intersection point.
Thus we use $\ell-1$ of the $x_j$'s in two edges and the last $x_j$ together with  $Z_{kr}$ and $(e'- e) \cap D$.
 Lastly, for $k \ge 5$ odd and $r=3$, we use a particular $Z_i$ twice to complete the odd cycle (since $|Z_i|=1$, this approach is valid only for $r=3$).

For $k\ge 5$ odd, suppose $u,v \in e \cap D$. Then again there is a path $P$ of length $k - 1$ in $H$ between $u$ and $v$
consisting of edges $Z_i \cup \{x_j,y\}$ with $y \in D$ such that $V(P) \cap e  = \{u,v\}$,
and $P \cup \{e\}$ is a $k$-cycle in $H$.

Finally, if $ k \ge 6$ is even, $r=3$,  $e=uvw,e'=u'v'w'$ with $e \cap e'= \emptyset$,
and $\{u,v,u',v' \} \subset D$, then we form a $C_k$ as follows: If $k=6$ we use the edges
$e, x_1z_1u, x_1z_2u', e', x_2z_3v', x_2z_4v$ where $Z_i=\{z_i\}$ for all $i$.
If $k>6$ then instead of the edge $x_2z_4v$, we use an edge $x_2z_4y$ for some $y\in D$, expand the path using the remaining $x_i$'s and $z_i$'s,
 and close the path with $x_{\ell}z_{2\ell}v$.  We obtain a cycle of length $2\ell+2=k$ as desired.
\qed

\medskip

{\bf Claim 4.} {\it $m > {n - 3r - 3k \choose r - 2}$.}

{\it Proof.} Suppose that $k$ is even and there are  $e,e' \in B$ with $|e \cap e'| = 1$.
Let $u \in e - e'$ and $v \in e' - e$ and let $f$ be an $r$-set with $f \cap (e \cup e') = \{u\}$ and $|f \cap L| = 1$.
If no such $r$-set is an edge of $H$, then $m \geq {n - |e \cup e' \cup L| \choose r - 2}$ and we are done. So we may assume that there is such an $f\in H$.
If $k>4$, then let $g$ be an $r$-set disjoint from $f$ and with $g \cap (e \cup e') = \{v\}$ and $|g \cap L| = 1$.
 If $k=4$, then let $g$ be an $r$-set with $g \cap (e \cup e' \cup f) = \{v\} \cup (f \cap L)$.
Let us argue that $g \not \in H$. Indeed, if $k>4$ and $g \in H$, then  we find a path $P$ of length $k - 2$  in $H$ as in Claim 3 containing $f$ and $g$,
and $P \cup \{e,e'\}$ is a $k$-cycle in $H$. If $k=4$, then $e, e', f, g$ is already a 4-cycle. Since $g \not\in H$ we have $g \in M$ and hence
$$m=|M| \ge {n - |e \cup e' \cup f \cup L|  \choose r - 2}>{n-3r-3k \choose r-2}.$$

If $r>3$, then by Frankl's theorem~\cite{Frankl}, $|B|>f(n,k,r)$ implies that there exist $e,e'\in B$ with $|e \cap e'|=1$. Now we are done by the preceding argument.
If $r=3$ and $k=4$, then by definition of $f(n,4,3)$ we find $e,e'$ with $|e \cap e'|=1$ and we are again done. If $r=3$ and $k \ge 6$ is even
and we cannot find such $e,e'$ with a singleton intersection, 
then  there are $e,e' \in B$ with $e\cap e'=\emptyset$ (this is easy to see since if we have more than $f(n, k, 3)=n-\ell-2$ triples on $n-\ell$ points and no singleton intersection, then we must have many disjoint complete 3-graphs on four points). Then for every $i$ and every $u \in e \cup e'$,  $d_H(x_iu)<3k$ for otherwise we can build a $k$-cycle using $e, e'$ and $k-2$ edges each containing some $x_i$ and at most one point 
of $e \cup e'$ (many of the edges will not intersect $e_1 \cup e_2$ if $k$ is large). 
This immediately gives at least $n-9-3k$ triples in $M$ that contain both $x_i$ and $u$ and Claim 4 is proved in this case.

If $k$ is odd,
then pick any edge $e \in B$ and apply a similar argument. \qed

For $0 \le i \le r$, define $B^r_i=\{e \in B: |e \cap (Y - D)|=i\}$.

\medskip

{\bf Claim 5.} {\it $|B^r_r|< \beta m$.}

{\it Proof.}  Recall that $c$ satisfies  $\ex_r(n, C_k)< cn^{r-1}$ for all $n \ge 1$. As $B^r_r$ itself has no $C_k$, we can apply this weaker bound to obtain
\begin{equation}\label{eqbr} |B^r_r| \le \ex_r(n - |D|,C_k) <
c(n - |D|)^{r-1}. \notag \end{equation}
Since $n$ is large, Claim 4 implies that $c2^{r-1}(\ell kr)^{r-1} < (\beta/2)m$ and
Claim 2 gives
$$|B^r_r| < c\left(\ell kr+ \frac{k^3r^rm}{n^{r-2}}\right)^{r-1}
< c2^{r-1}\left((\ell kr)^{r-1}+\left(\frac{k^3r^rm}{n^{r-2}}\right)^{r-1}\right)< \frac{\beta}{2} m+c'\frac{m^{r-1}}{n^{(r-2)(r-1)}}$$
where $c'=c2^{r-1}(r^rk^3)^{r-1}$. By (\ref{alpha})
and $m <\alpha n^{r-1}$,
$$c'\frac{m^{r-1}}{n^{(r-2)(r-1)}}=c'm\left(\frac{m}{n^{r-1}}\right)^{r-2}\le c'm\alpha^{r-2}<\frac{\beta}{2}m$$
and the claim follows. \qed
  \medskip

{\bf Claim 6.} {\it $|B^r_{r-1}|< \beta m$ for $r\ge 4$ and $|B^3_2|< 3m/4$.}

{\it Proof.}  Partition $B^r_{r-1}$ into $P^r \cup Q^r$, where $P^r$
comprises those $r$-sets $e \in B^r_{r-1}$ with  $d_{B^r_{r-1}}(e - D)=1$. Clearly $|P^r|<{|Y|-|D| \choose r-1}<(\beta/2) m$ as in Claim 5.

Let us now focus on $Q^r$.  Let $ F$ be the collection of $(r-1)$-sets $f\subset Y- D$ such that there exists $e \in B^r_{r-1}$ with $f \subset e$.
 We now partition the argument depending on whether $r=3$ or $r\ge 4$

Suppose that $r=3$. Then $ F$ is a (graph)  matching for if we have $vw$ and $vw'$ in $ F$,
then we have (by definition of $Q^3$) distinct vertices $y, y'$ and edges $vwy, vw'y'$ in $B^3_2$. This contradicts Claim 3. We will prove that $|Q^3|\le 2m/3$. Suppose for contradiction that $|Q^3|>2m/3$.
Then by averaging, there is a vertex $u \in D$
with $d_{B^3_2}(u) \ge \lceil 2m/(3n) \rceil: =t$. Let $v_1w_1, \ldots, v_{t}w_{t}$ be the neighbors of $u$ in $Q^3$
(meaning that $uv_iw_i \in Q^3$ for all $i$). Note that these pairs form a matching. Given $i<j$, there are at
least $2(|D|-2)$ sets of $M$ containing an element of $\{v_i, w_i\}$ or  at
least $2(|D|-2)$ edges of $M$ containing an element of $\{v_j, w_j\}$.  Indeed, if this is not the case, then we can form a copy
of $C_k$ using $uv_iw_i$ and $uv_jw_j$.  Since the pairs $\{v_iw_i\}_{i=1}^t$ form a matching this implies that $|M|\ge 2(|D|-2)(t-1)$. Since $m$ is large by Claim 4 and $\alpha$ is small this is at least $2 \times (0.9)n \times (\frac{2m}{3n}-1)>m$, contradiction.

Next suppose that $r\ge 4$. In this case $ F$ is a collection of $(r-1)$-sets on $D$ that have no singleton intersection by Claim 3.
  We conclude by a result of Keevash-Mubayi-Wilson~\cite{KMW} that $| F|<{n-|D| \choose r-3}$ and hence that
$$|Q^r| < | F| n < {n-|D| \choose r-3} n.$$
By Claim 2, there exists  $C$ depending
only on $k$ and $r$ such that this is at most
$$C n \left(\frac{m}{n^{r-2}}\right)^{r-3}=Cm \frac{m^{r-4}}{n^{(r-2)(r-3)-1}}.$$
Since $m < n^{r-1}$, $(r-1)(r-4)<(r-2)(r-3)-1$ and $n$ is large, the last expression is at most $(\beta/2) m$ and the claim follows. \qed

Since $|B|=m+f(n, k, r)$, Claims 5 and 6 imply that
$|B^r_{r-1}|+|B^r_r| < (2\beta+3/4) m < m$ and therefore
 $|B^r_0 \cup \ldots \cup B^r_{r-2}|> f(n,k,r)$.

If $k$ is odd, then $B^r_0 \cup \ldots \cup B^r_{r-2} \ne \emptyset$.
If $k$ is even and  $r \ge 4$ then there are edges $e, e' \in B^r_0
\cup \ldots \cup B^r_{r-2}$ such that $|e \cap e'| = 1$. This is because for $r \ge 4$ the extremal function
for $P_2$ is the same as the extremal function for $\{P_2, 2P_1\}$ by~\cite{Frankl} as long as $n$ is sufficiently large
(in both cases the extremal example is obtained by taking all $r$-sets that intersect a specific set of two points).
If $(k,r)=(4,3)$, then  by definition of $f(n,4,3)$ there are edges $e, e' \in B^3_0
\cup B^3_{1}$ such that $|e \cap e'| = 1$.
Finally, if $k \ge 6$ is even, $r=3$ and $|B^3_0 \cup  B^3_{1}|> f(n,k,3)=n-\ell-2$ then we  find two edges $e,e' \in B^3_0 \cup B^3_1$ with $|e \cap e'|\le 1$.
In all four cases above we contradict Claim 3. This completes the proof of Theorem \ref{main}. \qed

\subsection{Part IIIb : Exact result for paths.\mbox{ }}\label{exactpaths}
We closely follow the proof  in Section~\ref{exact} except that we  replace $f(n,k,r)$ by $h(n,k,r)$,
where $h(n,k,r) = 0$ if $k$ is odd and $h(n,k,r)=ex_r(n-\ell, \{P_2, 2P_1\})$ if $k$ is even.
   Claims 1, 2 and 5
follow immediately and Claim 4 follows by a very similar proof.  We strengthen Claim 3 as follows.

\medskip

{\bf Claim 3$'$.} {\it No two edges $e,e' \in B$ have $|e \cap e'|\le 1$, $(e -e') \cap D \neq \emptyset$ and
$(e' - e) \cap D \neq \emptyset$.
If $k$ is odd, then no edge $e \in B$ has $|e \cap D| \geq 1$.}

{\it Proof.} In the first case, we may form a path using the two vertices of $e \triangle e'$ in $D$ and $2\ell$ other edges.
This is a  path of length $2\ell+2\ge k$. In the case when $k$ is odd, we form a path of length $2\ell+1=k$ ending at $e$ by the same procedure. \qed

\medskip

If $k$ is odd, then Claim 3$'$ implies that $B=B^r_r$ and Claim 5 implies the contradiction  $m\le |B|<\beta m$. Let us suppose that $k$ is even.
We now observe that Claim 6 also holds (in fact we can improve the argument when $r=3$ to obtain $4(|D|-1)$ instead of $2(|D|-1)$ as it is easier to form  a $k$-path), so
$|B^r_0 \cup \ldots \cup B^r_{r-2}|>h(n,k,r)$ and we find a $P_2$ or a $2P_1$ in this union. This contradicts Claim 3$'$ and completes the proof. \qed

\section{Proof of Theorem \ref{loosemain}}\label{proofofloosemain}

In this short section we show how to modify the proof of Theorem~\ref{main} to prove Theorem~\ref{loosemain}.
The case of minimal paths is easier than minimal cycles, so we concentrate only on minimal cycles.
We only prove the case $r=3$ as all other cases are covered by the result of F\"uredi-Jiang~\cite{FJ} (though our proof works just as easily for all $r\ge 3$ and $k \ge 4$).
We closely follow the proof of Theorem \ref{main}.  We may assume that $k \ge 4$ is even as the case $k=3$ is already solved
in~\cite{CK, FF, MV} and if $k \ge 5$ is odd, then we apply Theorem \ref{main} directly. Since
$C_k \in \CC_k$, we immediately obtain a stability result
(Theorem \ref{stab}) for $\CC_k$. Now we repeat the proof in Section \ref{exact} with $f(n,k,r)$ replaced by $f(k)$,
where $f(k)=0$ if $k$ is odd, $f(k)=\lfloor (n-1)/r\rfloor$ if $k=4$ and $f(k)=1$ if $k\ge 6$ is even.
The proofs of Claims 1, 2, 4, 5 and 6 remain the same or very similar and we do not repeat them.  Claim 3 can be strengthened by replacing $|e \cap e'|=1$ with
$|e \cap e'|\ge 1$ since it is enough to find a minimal cycle.

Suppose that $k=4$, $\ell=1$ and we are trying to find a minimal 4-cycle. Then $|B^3_2| + |B^3_3|<
 (\beta+3/4)m\le  (1/10+3/4)m < (6/7)m$ and therefore
 $|B^3_0|+|B^3_1|=m+f(k)-|B^3_2|-|B^3_3|>f(k)+m/7$.  If $|B^3_0|>f(k)$, then we find
 $e, e' \in B^3_0$ with $e \cap e' \ne\emptyset$ which contradicts (the strengthened) Claim 3.  So we may assume that $|B^3_0|\le f(k)$ and $|B^3_1|>m/7$.  Each edge of $B^3_1$ has a
 vertex in $Y- D$, and since $n$ is large,
 $|Y- D|< m/7$.  Therefore
 there is a vertex $v \in Y-D$ with $d_{B^3_1}(v)>1$. This again contradicts Claim 3.

 Now we suppose that $k \ge 6$ is even, and $f(k)=1$.
 If $|B^3_0|>f(k)=1$, then there are two  edges $e,e' \subset D$ and this contradicts Claim 3 (no matter what their intersection size).
 We may therefore assume  that $|B^3_1|>m/7$ and this again contradicts Claim 3 as above. \qed

\bigskip

{\bf Acknowledgment.} We thank Zoltan F\" uredi and Tao Jiang for  helpful comments.





\end{document}